\documentclass[11pt]{article}
\usepackage{a4,amssymb,verbatim,amsmath}

\newtheorem{theorem}{Theorem}[section]

\newtheorem{remark}[theorem]{Remark}

\newtheorem{ex}{Example}[section]

\newtheorem{ass}{Assumption}[section]

\numberwithin{equation}{section}

\begin{document}

\def\sso#1{\ensuremath{\mathfrak{#1}}}  
\def\ssobig#1{{{#1}}}

\newcommand{\ddx}{\partial \over \partial x}
\newcommand{\ddy}{\partial \over \partial y}

\begin{center}
{\bf \Large Conservative difference schemes \\ for one-dimensional flows of polytropic gas}
\end{center}

\bigskip

\begin{center}
{Roman Kozlov}

\bigskip
{Department of Business and Management Science, \\
Norwegian School of Economics, \\ 
Helleveien 30, 5045, Bergen, Norway;  \\ 
{e-mail: Roman.Kozlov@nhh.no}}
\end{center}


\bigskip


\begin{center}
{\bf Abstract}
\end{center}
\begin{quotation}
The paper considers one-dimensional flows of polytropic (calorically ideal) gas. 
These flows include three cases of gas dynamics: 
plain one-dimensional flows (one-dimensional space),  
radially  symmetric flows in two-dimensional space 
and   spherically  symmetric  flows  in three-dimensional space.  
Starting with the difference schemes which have conservation laws of mass and energy 
(as well as conservation of momentum and the center of mass motion for the  plain one-dimensional flows), 
we find difference schemes  which also have additional conservation laws 
for the special values of the adiabatic exponent $\gamma = 1 + 1 /d $, where $d$ is the space dimension.  
\end{quotation}

\bigskip

{\bf Key words:} 

Conservative difference schemes

Polytropic gas

Mass Lagrange coordinate


\section{Introduction}      

\label{Introduction}

In many applications gas phenomena can be modeled using the ideal gas 
which is characterized by the equation of state~\cite{Chernyi,  Landau, Ovsiannikov,  Samarskii, Toro}   
\begin{equation}      \label{idea1}
p = \rho R T  . 
\end{equation}   
Here  $ p$ is the pressure, $ \rho $ is the density and  $ T $ is the pressure. 
Constant $R$   is called the specific gas constant.

We will assume that the specific internal energy  $ \varepsilon  $ 
 is a linear function of the temperature (for constant volume gas processes)  
\begin{equation}      \label{idea2}
\varepsilon (T) = C_v T = { RT \over \gamma - 1}  . 
\end{equation}   
Here   $ \gamma $  is the adiabatic exponent.   
It is defined the ratio of specific  heat capacities  
\begin{equation*}   
\gamma = { C_p \over C_v } ,  
\end{equation*}   
where  $C_p$ is the specific heat capacity at constant pressure 
and   $C_v$ is the specific heat capacity at constant volume.


The ideal gas~(\ref{idea1}) with the assumption~(\ref{idea2}) is called {\it polytropic}   
(also known as {\it calorically ideal}). 
Excluding temperature, we obtain 
\begin{equation}        \label{ideal_gas}
\varepsilon   = { 1 \over \gamma  - 1 }   { p \over \rho}    .  
\end{equation}
It is well known~\cite{Landau, Toro} that for polytropic gases specific heat capacities 
$C_v$ and $ C_p$,  and consequently  the adiabatic exponent $ \gamma $ 
are constant.

In this paper we treat one-dimensional flows 
for the  polytropic  gas~(\ref{ideal_gas}). 
For physical applications it makes sense to restrict to the case $\gamma > 1$.  
However, the results of the paper are valid for  $ \gamma  \neq 0, 1   $.
We consider three types of the gas flows, 
namely,  plain one-dimensional flows (one-dimensional space),  
radially symmetric flows  in two-dimensional space 
and spherically  symmetric  flows  in three-dimensional space.  
It is convenient to examine these three cases together. 
We will call them {\it one-dimensional}  flows.


Equations of the gas dynamics are derived from 
the conservation of mass,  momentum and energy~\cite{Chorin, Samarskii}.  
Therefore it is reasonable to preserve these conservation laws in difference schemes.   
The purpose of the paper is to find difference schemes 
which preserve the  conservation laws of the gas dynamics equations.  
Not only the basic conservation laws of mass and energy 
(as well as momentum and motion of the center of  mass for  the  plain one-dimensional flows),  
but also additional conservation laws which exist for the special values of $\gamma = 1 + 1 /d $, 
where $d$ is the space dimension.   
To achieve this goal we will use the mass Lagrange coordinates.  
For  the plain one-dimensional flows such difference scheme on a uniform spatial mesh was 
suggested  in~\cite{Korobitsyn}.  
This difference scheme will be recovered as a particular case.   
We stress that the schemes to be constructed preserve 
the conservation laws which hold for the weak solutions 
(i.e., solutions with discontinuities) of the gas dynamics equations.  
The conservation law which  hold only for the smooth solutions, namely~(\ref{Euler_state}), will be discarded.

The paper belongs to a fairly new research direction in numerical methods which was called  
{\it structure-preserving numerical methods}  
(also known as  {\it geometric numerical integration}). 
Though the field is mainly concern with the qualitative properties of ODES~\cite{Hairer}, 
there are also applications to PDEs~\cite{Furihata, Shashkov}.

Numerical methods for gas dynamics include many different approaches. 
The methods differ in the choice of coordinates, the way to treat discontinuous solutions, 
approximation order, stability properties, etc. 
For review of the methods we refer to books~\cite{LeVeque, Laney,  Toro} and references wherein. 
The schemes developed in the present paper can be described as uniform schemes  
(no isolation and tracking of shocks) in Lagrangian coordinates.

The paper is organized as follows:  
In the next section 
we describe the gas dynamics equations of the 
one-dimensional flows in Euler coordinates and their conservation laws.  
We introduce the mass Lagrangian coordinate 
and rewrite the gas dynamics equations and their conservation  laws in Lagrange coordinates in section~\ref{Lagrange}.  
Then, in section~\ref{scheme}  we obtain the conservative difference schemes. 
The results of the paper are summed up in the final section~\ref{conclusion}.


\section{One-dimensional gas flows in Euler coordinates}      

\label{Euler}

\subsection{Equations of one-dimensional flows} 

\label{Euler_system}

In Euler coordinates the gas dynamics equations of the one-dimensional  flows 
can be presented as~\cite{Samarskii, Toro}  
\begin{subequations}     \label{GD_Euler_coordinates}
\begin{gather} 
 \rho   _t   +   u   \rho _r   +  {  \rho   \over r^n  }    ( r^n u )  _r     =   0    ,    \label{Euler_rho}
\\
 u   _t   +   u   u _r   +  {  1 \over  \rho   }    p_r     =   0     ,      \label{Euler_u}
\\
   \varepsilon    _t   
+ u    \varepsilon    _r 
+  {  p  \over \rho  r^n  }    ( r^n u  )  _r       =   0     .       \label{Euler_e}
\end{gather} 
\end{subequations} 
Here  we distinguish the case $n = 0$ with coordinate  $- \infty <r < \infty$  and velocity $u$ 
from the cases $n = 1,  2 $  with radial distance from the origin $0< r < \infty$   and the radial velocity $u$.      
The system should be supplemented by the equation of state 
\begin{equation}       \label{energy}
\varepsilon = \varepsilon ( \rho , p )  . 
\end{equation}   
We have $n = 0, 1, 2$ 
for  the plain one-dimensional flows, 
the radially symmetric two-dimensional flows 
and the spherically  symmetric three-dimensional flows, respectively.  
Note that for these cases  $ n = d - 1$, where $d= 1, 2, 3$ is the space dimension.

For the  polytropic gas~(\ref{ideal_gas}),  
which we consider in this paper, 
we can replace the last equation~(\ref{Euler_e})  with 
\begin{equation}     \label{Euler_p}
 p _t   +   u   p _r   +  {  \gamma  p   \over r^n  }    ( r^n u )  _r     =   0    . 
\end{equation}
From the equations~(\ref{Euler_rho})  and~(\ref{Euler_p})  
it follows that 
\begin{equation}     \label{state_CL_Euler}
\left( { p \over\rho  ^{\gamma} }  \right)_t 
+ u  \left( { p \over\rho  ^{\gamma} }  \right)_r = 0   . 
\end{equation}  
It leads to the conservation of the entropy.   
For the polytropic  gas~(\ref{ideal_gas}) the entropy is given as (see~\cite{Chernyi,  Chorin, Landau, Ovsiannikov, Samarskii})
\begin{equation}     \label{Euler_entropy_2}
S = { R \over \gamma -1 }  \ln \left( { p \over\rho  ^{\gamma} } \right)       . 
\end{equation}  
Its conservation along the streamlines of the gas flow follows from~(\ref{state_CL_Euler}), 
which can be rewritten as 
\begin{equation}    \label{Euler_entropy_1}
S_t + u S_r = 0 .  
\end{equation}



\subsection{Conservation laws}      

\label{Euler_CL}

The conservation laws  of the system~(\ref{GD_Euler_coordinates}) for the gas~(\ref{ideal_gas}) 
can be obtained by direct computation. They can be split into several groups as follows.

\begin{enumerate}

\item 

General case

In the general case we get two conservation laws:

\begin{itemize}

\item   

Conservation of mass
\begin{equation}      \label{Euler_mass}
\left[ r ^n \rho  \right]   _t   +     \left[ r^n  \rho  u  \right]   _r  = 0   ; 
\end{equation}

\item   

Conservation of energy 
\begin{equation}      \label{Euler_energy}
\left[ r  ^n   \left(  \rho   \varepsilon +  { \rho u^2 \over 2 }    \right)  \right]   _t   
+   \left[ r^n   \left(  \rho   \varepsilon +  { \rho u^2 \over 2 }    + p   \right)  u  \right]    _r  = 0   .
\end{equation}

\end{itemize}

\item

Plain one-dimensional flows

For  $ n= 0  $ there are two more  conservation laws:

\begin{itemize}

\item 

Momentum 
\begin{equation}        \label{Euler_momentum}
\left[
   \rho    u  
\right]   _t   
+   \left[ 
   \rho   u  ^2  +   p 
 \right]    _r  = 0    ; 
\end{equation}

\item 

Motion of the center of mass 
\begin{equation}      \label{Euler_center}
\left[
   \rho   ( r - t  u  ) 
\right]   _t   
+   \left[ 
   \rho   u  ( r  - t u ) - t p 
 \right]    _r  = 0     .
\end{equation}

\end{itemize} 
These two conservation laws hold only for the plain one-dimensional flows. 
For the higher dimensional cases of gas dynamics 
with a radial or spherical symmetry 
the conservation of momentum  and no motion of the center of mass 
are built into the model differential equations, 
which are obtained from two- and three-dimensional gas dynamics equations using the symmetry reduction procedure.

\item 

Special values of $\gamma$

For $\gamma = { 1  + 2 / d } $, 
i.e. for  $ \gamma =  3$, $2$ and $ 5/3$ in the cases of dimension $d = 1$, $2$ and $3$, 
there are two additional  conservation laws 
\begin{multline}       \label{Euler_additional_1}
\left[      r^n \left(   2t    \left(  \rho   \varepsilon +  { \rho u^2 \over 2 }    \right)    
-   r    \rho u     \right)   \right]   _t   
\\
+   \left[    r^n   \left( 2t    \left(  \rho   \varepsilon +  { \rho u^2 \over 2 }    + p   \right)  u 
- r     ( \rho u^2 + p   )  
 \right) \right]    _r  = 0    
\end{multline} 
and 
\begin{multline}      \label{Euler_additional_2}
\left[   r  ^n      \left(    t^2      \left(  \rho   \varepsilon +  { \rho u^2 \over 2 }    \right)    
-  t  r     \rho u  
+ { r^2  \over 2 }      \rho 
  \right)   \right]   _t   
\\
+   \left[ r^n   \left(    t^2     \left(  \rho   \varepsilon +  { \rho u^2 \over 2 }    + p   \right)  u 
- t r    ( \rho u^2 + p   )  
+ { r^2  \over 2 }     \rho   u   
 \right)  \right]    _r  = 0     .
\end{multline} 
They do not have a clear physical interpretation as the other conservation laws.

\end{enumerate}

Let us note that the conservation laws~(\ref{Euler_mass})-(\ref{Euler_center}) 
hold for equations~(\ref{GD_Euler_coordinates}) with any equation of state~(\ref{energy}). 
The additional conservation laws~(\ref{Euler_additional_1}) and~(\ref{Euler_additional_2})
hold only  for the  polytropic   gas~(\ref{ideal_gas}) 
with the special values $\gamma$.

\begin{remark} 
For the gas~(\ref{ideal_gas}) 
there is also  the conservation law
\begin{equation}       \label{Euler_state}
\left[ r ^n \rho   F \left( { p   \over \rho ^{ \gamma } } \right) \right]   _t   
+     \left[ r^n  \rho  u    F \left( { p   \over \rho ^{ \gamma } } \right) \right]   _r  = 0  , 
\end{equation}
where $F$ is a differentiable function. 
From 
\begin{multline*}      
\left[ r ^n \rho   F \left( { p   \over \rho ^{ \gamma } } \right) \right]   _t   
+     \left[ r^n  \rho  u    F \left( { p   \over \rho ^{ \gamma } } \right) \right]   _r  
\\
= F \left( { p   \over \rho ^{ \gamma } } \right)  
\left( 
\left[ r ^n \rho    \right]   _t   
+     \left[ r^n  \rho  u    \right]   _r  
\right)   
+ 
r ^n \rho   
F ' \left( { p   \over \rho ^{ \gamma } } \right) 
\left[  
\left(   { p   \over \rho ^{ \gamma } }  \right)    _t   
+   u   \left(    { p   \over \rho ^{ \gamma } } \right)    _r 
\right] 
\end{multline*}      
we see that it  holds due to the mass conservation~(\ref{Euler_mass})  
and the entropy conservation property~(\ref{Euler_entropy_1}).

The entropy is conserved along the streamlines only for smooth solutions. 
It is not conserved for discontinuous solutions (shocks)~\cite{Chernyi,  Chorin, Ovsiannikov}.  
Since in this paper we concern with conservation laws which hold for the weak solutions 
we will not consider property~(\ref{Euler_state}) important 
to preserve under discretization. 
\end{remark}


We presented all conservation laws of the one-dimensional flows of the gas dynamics equations. 
One can find them by direct computation 
or   by an appropriate reduction of the three-dimensional conservation laws. 
Conservation laws of three-dimensional gas dynamics were obtained  
in~\cite{Ibragimov2} (see also~\cite{Ibragimov}) 
with the help of a variational formulation and Noether's theorem 
(it requires some assumptions) 
and by direct computation in~\cite{Terentev}.    
Among the 13  conservation laws of the three-dimensional case 
all but one can be integrated over discontinuities~\cite{Terentev}.  
The only  conservation law which cannot be integrated over discontinuities gets reduced to~(\ref{Euler_state}) 
in the case of one-dimensional flows. 
It cannot be integrated over discontinuities because the entropy is  not 
conserved for the discontinuous solutions~\cite{Chorin, Landau}.    
In~\cite{Ibragimov2, Ibragimov} one can find a symmetry interpretation of the conservation laws, 
i.e.  the correspondence between the conservation laws and Lie point symmetries 
of the three-dimensional gas dynamics equations.

In this paper we will find difference schemes 
with difference  analogs of the conservation laws~(\ref{Euler_mass})-(\ref{Euler_additional_2}). 
To achieve this goal we introduce Lagrange coordinates.


\section{One-dimensional gas flows in Lagrange coordinates}      

\label{Lagrange}

For preservation of the conservation laws under discretization 
we will make use of the Lagrange coordinates. 
Let us introduce the mass Lagrange variable as 
\begin{equation}      \label{Lagrange_coordinate}
s = \int  _{r_0}    ^r   \rho( t , y)  y^n   dy   
\end{equation}
and rewrite the gas dynamics equations and their conservation laws in the Lagrange coordinates.

The total derivatives  $D_t ^L   $ and $ D_s $ 
with respect to $t$  and $s$ in  the Lagrange coordinates  are related 
to the total derivatives  $D_t ^E  $ and $ D_s $ 
with respect to time $t$ and $r$ in the Euler coordinates as 
\begin{equation} 
D_t ^L  
= D_t  ^E  + u   D_r , 
\qquad 
D_s = { 1 \over r^n \rho}   D_r   . 
\end{equation}
We remark  that Lagrange time derivative  $ D_t ^L  $  
represents the material time derivative, 
which corresponds to the differentiation along the streamlines.

\subsection{Equations of one-dimensional flows} 

\label{Lagrange_system}

The gas dynamics equations~(\ref{GD_Euler_coordinates})  
are transformed from  the Euler coordinates $ (t,r)$ to the Lagrange coordinates $ ( t,s) $ as 
\begin{subequations}     \label{GD_Lagrange_equations}
\begin{gather}
   \rho   _t   +    \rho ^2    ( r^n  u )_s   = 0   ,     \label{Lagrange_rho}
  \\ 
u _t   +    r^n   p_s    =  0  ,   \label{Lagrange_u}
   \\
\varepsilon  _t   +   p (   r^n   u  ) _s   = 0   .   \label{Lagrange_e}
\end{gather}
\end{subequations}
The Euler spatial coordinate $r$ is given in the Lagrange coordinates by the equations
\begin{subequations}     \label{Lagrange_xt_xs}
\begin{gather}
r_t =    u   ,       \label{Lagrange_xt}
  \\
r_s =  { 1 \over  r^n \rho}        \label{Lagrange_xs}   .
\end{gather}
\end{subequations}

In case of the polytropic ideal gas~(\ref{ideal_gas}) 
the last equation~(\ref{Lagrange_e}) can be replaced by the equation for the pressure 
\begin{equation}    \label{Lagrange_p}
p _t  +  \gamma  \rho p  ( r^n  u )_s  =  0    .
\end{equation}
It is easy to see in the Lagrange case  
that the conservation of the entropy along the streamlines follows from 
\begin{equation}     \label{state_CL_Lagrange}
\left( { p \over\rho  ^{\gamma} }  \right)_t  = 0   , 
\end{equation}  
a consequence of  the equations~(\ref{Lagrange_rho}) and~(\ref{Lagrange_p}).


From equations~(\ref{Lagrange_rho}) and~(\ref{Lagrange_e}) 
it follows that 
\begin{equation}    \label{work}
\varepsilon  _t =   -  p \left(  { 1 \over \rho }  \right) _t      . 
\end{equation}
This relation shows that the change of the specific internal energy   $ \varepsilon $  
is caused by the work of the pressure forces. 
It is important to preserve this relation 
in addition to the balance of the total energy 
under discretization   
for qualitatively correct numerical simulation~\cite{Samarskii}. 
It keeps the correct balance between the  specific internal energy and the specific kinetic energy.


\subsection{Conservation laws}

\label{Lagrange_CL}

Let us rewrite the conservation laws of point~\ref{Euler_CL} 
for  the Lagrangian coordinates.

\begin{enumerate}

\item 

General case

In the general case there are  two  conservation laws:

\begin{itemize}

\item 

Conservation of mass  
\begin{equation}      \label{Lagrange_mass}
\left[  { 1 \over \rho }  \right]_t 
- 
 [ r^n  u ]_s   
= 0   ; 
\end{equation}

\item 

Conservation of energy 
\begin{equation}      \label{Lagrange_energy}
\left[ \varepsilon   + {1 \over 2} u^2 \right]   _t  
+ 
   [  r^n p   u  ] _s   
=   0    .
\end{equation}

\end{itemize}

\item

Plane    one-dimensional flow

For  $ n= 0  $  there are also

\begin{itemize}

\item 

Conservation of momentum 
\begin{equation}      \label{Lagrange_momentum}
\left[    u \right]    _t  
+ 
  [   p ]  _s  
=   0     ; 
\end{equation}

\item 

Motion of the center of mass 
\begin{equation}     \label{Lagrange_center}
\left[ r  -   t  u       \right]    _t  
 -    [    t p ]  _s  
=   0    .
\end{equation}

\end{itemize}

\item

Special values of $\gamma$

For the ideal gas~(\ref{ideal_gas})  with $ \gamma =  { 1 + 2  /  d } $  
there are two additional conservation laws  
\begin{equation}       \label{Lagrange_additional_1}
\left[ 2t  \left( \varepsilon   + {1 \over 2} u^2 \right)   - r u \right]    _t  
+ 
  [     r^n    p     ( 2t    u   - r )    ]  _s  
=   0 
\end{equation}
and 
\begin{equation}      \label{Lagrange_additional_2}
\left[ t^2   \left( \varepsilon   + {1 \over 2} u^2 \right)   -   t r u    +  { r^2 \over 2 }   \right]    _t  
+ 
  [    r^n p (     t^2     u   - t r    ) ]  _s      
=   0   . 
\end{equation}

\end{enumerate}

\begin{remark} 
For the gas~(\ref{ideal_gas}) we also have  the conservation law
\begin{equation}       \label{Lagrange_state}
\left[ { p   \over \rho ^{ \gamma } } \right]   _t    = 0  , 
\end{equation}
which is the analog of~(\ref{Euler_state}). 
\end{remark}


\section{Conservative difference schemes}      

\label{scheme}

\subsection{Notations}

\label{scheme_notations}

For the mass Lagrange variable we introduce a mesh with points  $ s_i $, $ i = 0, 1, ..., N$.  
Generally, this mesh is nonuniform. We denote the mesh  {step}lengths    as 
\begin{equation*}   
h_i = s_{i+1} - s _i . 
\end{equation*}   
Since the difference schemes will consists of equations for two time layers, 
namely $t_ {j} $ and $t_ {j+1} $, 
we can denote the  time {step}length as  $ \tau  $. 
Of course, we can take  nonuniform {step}lengths  $  \tau _j $, 
which generate the time mesh points $ t_ {j+1} = t_j  +  \tau _j $, $ j = 0, 1, ...$

The kinematic variables $r $ and $u$  are taken 
in the nodes $(t_j  , s_i)$ of the two-dimensional mesh. 
We denote them as 
\begin{equation*}   
u_- =  u_{i-1} ^j  , 
\quad
u =  u_ i  ^j , 
\quad  
u_+ =  u_{i+1} ^j  , 
\quad  
\hat{u} _- =  u_{i-1} ^{j+1}   , 
\quad
\hat{u}  =  u_ i  ^{j+1}  , 
\quad  
\hat{u} _+ =  u_{i+1} ^{j+1}    . 
\end{equation*}   
Their difference   time derivative and forward and backward difference spatial derivatives will be  
\begin{equation*}   
u_t = { \hat{u}   -  u   \over  \tau  } , 
\qquad
u_s = { u_{i+1}  ^j   - u _i ^j    \over  s_{i+1} - s_i   }    
=  { u_{+}   - u   \over   h_i   }  , 
\qquad
u_{\bar{s}} = { u _i ^j    - u_{i-1}  ^j  \over  s_i   - s _{i-1} }  
=  { u   - u_-   \over   h_{i-1}   }   .
\end{equation*}   
It is helpful to use a special notation for the average value of two function values 
taken in  the neighbouring nodes of the same time layer 
\begin{equation*}   
< f (u,r) > = { f (u,r)  + f (u_+,r_+)  \over 2}   .
\end{equation*}

The thermodynamical  variables $\rho  $, $p $  and $ \varepsilon $  
are assigned to the midpoints $(t_j , s_{ i+1/2 }) $, $  s_{ i+1/2 } = (  s_{ i  } +  s_{ i+1 } ) / 2 $.   
For example, 
\begin{multline*}   
\rho_-   =  \rho  _ {i-1/2}  ^j  , 
\quad  
\rho   =  \rho  _ {i+1/2}  ^j  , 
\quad  
\rho  _+ =  \rho  _ {i+3/2} ^j  , 
\\
\hat{\rho}_-   =  \rho  _ {i-1/2}  ^{j+1}  , 
\quad  
\hat{\rho}   =  \rho  _ {i+1/2}  ^{j+1}   , 
\quad  
\hat{\rho}  _+ =  \rho  _ {i+3/2} ^{j+1}    . 
\end{multline*}      
For their forward and backward difference spatial derivatives 
we take into account that the mesh is not uniform
\begin{equation*} 
 p_{\bar{s}} = { p _{ i+1/2 }  ^j -   p _{ i-1/2 }  ^j \over    { 1 \over 2} ( h_i + h_{i-1}  )   } , 
\qquad 
p_{s} = { p _{ i+3/2 } ^j  -   p _{ i+1/2 }  ^j \over    { 1 \over 2} ( h_{i+1} + h_{i}  )   }  .  
\end{equation*}
We will need the linear interpolation value of the pressure in the nodes. 
For example,  for the node $ ( t_j ,   s_i ) $ we get 
\begin{equation*} 
   p_* = ( {p_*} ) _{i}  ^j   = {   h_i  p _{ i-1/2 }  ^j   +  h _{i-1}    p _{ i+1/2 }  ^j   \over   h_i + h_{i-1}     }    .
\end{equation*}

For all variables we denote weighted values for the two neighbouring  time layers  as  
\begin{equation*} 
 y    ^{(\alpha)}    =  \alpha   \hat{y} + ( 1 - \alpha ) y   , 
\qquad 
0 \leq \alpha \leq 1   . 
\end{equation*}


\subsection{Discrete scheme} 

\label{scheme_system}

As a starting point for discretization of  
the equations~(\ref{GD_Lagrange_equations}) and~(\ref{Lagrange_xt}) 
we take the scheme~\cite{Samarskii}  
\begin{subequations}       \label{equation_preliminary}
\begin{gather}
\left(  { 1 \over \rho }  \right) _t = ( R  u  ^{(0.5)}   )_s      ,   \label{equation_discrete_mass}
\\
u _t =   -   R    p_{\bar{s}}    ^{(\alpha)}   ,      \label{equation_discrete_velocity}
\\ 
\varepsilon  _t =   -  p ^{(\alpha)} (   R    u  ^{(0.5)} ) _s    ,     \label{equation_discrete_energy}
\\
r_t =    u  ^{(0.5)}     .     \label{equation_discrete_coordinate}
\end{gather}
\end{subequations}
Here 
\begin{equation*} 
R = { \hat{r} ^{n+1} -    {r} ^{n+1}  \over ( n+1 ) ( \hat{r}  - r )  } 
= 
\left\{ 
\begin{array}{ll} 
1 , &   n  =  0 ; \\
 &  \\
{\displaystyle   { \hat{r} +  r  \over 2 }   }  , &   n  =  1 ; \\ 
 &  \\
 {\displaystyle   { \hat{r}^2  +   \hat{r} r  +  r^2  \over 3 }   }      , &   n  =  2 \\ 
\end{array} 
\right. 
\end{equation*}
is taken for  a weighted discretization of  $  r^n $.  
Note that    $ R \rightarrow r^n $ as  $ \tau \rightarrow 0$.  
If the spatial mesh is uniform with {step}length $h$, 
this difference scheme approximates the underlying differential equations with $ O ( \tau + h^2 )$.

The scheme~(\ref{equation_preliminary}) 
has four equations for five variables:    $\rho$, $u$, $\varepsilon$, $r$ and  $p$.  
It should be supplemented by a discrete equation of state, 
a discrete analog of~(\ref{energy}).  
For example, it can be taken in the same form that means 
\begin{equation*}
\varepsilon   _{ i+1/2}  ^j  =  \varepsilon     ( \rho   _{ i+1/2}  ^j ,  p  _{ i+1/2}  ^j   )     . 
\end{equation*}
However, we will not impose any discrete equation of state at the moment. 
The freedom to choose a discretization of the equation of state 
will be used to impose one additional conservation law.


We stress that the scheme has qualitatively correct discretization of the relation~(\ref{work}), 
namely 
\begin{equation}     \label{discrete_work}
\varepsilon  _t =   -  p ^{(\alpha)}   \left(  { 1 \over \rho }  \right) _t     .
\end{equation}


\begin{remark}
To derive the difference scheme~(\ref{equation_preliminary})  
one can start with a nonconservative discretization of 
the equations~(\ref{GD_Lagrange_equations}) and~(\ref{Lagrange_xt})  
with many free parameters.

In~\cite{Popov} 
(slightly  different approach  was undertaken in~\cite{Moskalkov1, Moskalkov2})   
the authors considered the case  $n=0$. 
Here one can start with the nonconservative discretization 
\begin{subequations}         \label{original}
\begin{gather}
\left(  { 1 \over \rho }  \right) _t =  u  ^{( \sigma _1 )}   _s      ,       \label{original_rho}
\\
u _t =   -       p_{\bar{s}}    ^{( \sigma _2 )}   ,        \label{original_u}
\\
\varepsilon  _t =   -  p ^{( \sigma _3 )}     u  ^{( \sigma _4 )}  _s    ,      \label{original_e}
\\
r_t =    u  ^{( \sigma _5 )}     .       \label{original_r}
\end{gather}
\end{subequations}
They imposed conservation laws of mass and energy that left only two free parameters.  
Note that conservation of momentum is provided by the equation~(\ref{original_u}).  
In~\cite{Samarskii}  they also required the scheme to keep the relation~(\ref{work}) 
under discretization and specified the scheme as given by~(\ref{equation_preliminary}) for $n=0$.

The scheme~(\ref{equation_preliminary}) for $n=0$ 
was extended to the other cases of the one-dimensional flows, 
i.e. for $ n =2,3$, with the help of a suitable discretization  of $r^n$ in~\cite{Samarskii}. 
\end{remark}



\subsection{Discrete conservation laws and discrete equation of state} 

\label{scheme_CL}

First, we provide the conservation laws of the difference scheme~(\ref{equation_preliminary}) 
which hold independently of the equation of state choice.

\begin{enumerate}

\item 

The general case

In the general case the scheme has two conservation laws:

\begin{itemize}

\item 
Conservation of mass is given by the equation~(\ref{equation_discrete_mass}).

\item

Conservation of energy 
\begin{equation}         \label{scheme_energy} 
\left[ \varepsilon   + {  < u^2 >   \over 2}       \right]   _t  
+ 
   [ R   p _{*} ^{(\alpha)}    u ^{(0.5)}  ] _s   
=   0  . 
\end{equation}

\end{itemize}

\item

Plane one-dimensional flows

For $n=0$ there are two more  conservation laws:

\begin{itemize}

\item 

Conservation of momentum   
\begin{equation}        \label{scheme_momentum} 
\left[ 
    u      
 \right]    _t  
+    
  [      p  ^{(\alpha)}    ]  _{\bar{s}}  
=   0     ; 
\end{equation}

\item 

Motion of the center  of mass 
\begin{equation}        \label{scheme_center} 
\left[ 
r  -   t   u      
 \right]    _t  
-  
  [    t   ^{(0.5)}    p  ^{(\alpha)}    ]  _{\bar{s}}  
=   0    , 
\end{equation} 
which was not mentioned by the authors of~\cite{Popov, Samarskii}.

\end{itemize}

\end{enumerate}

We can formulate the properties of the scheme as the following theorem.

\begin{theorem}    \label{first_theorem} 
The scheme~(\ref{equation_preliminary}) with an arbitrary equation of state~(\ref{energy}) 
possesses the conservation laws of mass~(\ref{equation_discrete_mass}) 
and energy~(\ref{scheme_energy}). 
It satisfies the discrete relation~(\ref{discrete_work}), 
which insures the correct balance of the specific internal and specific kinetic energies. 
For $n=0$ there are also conservation laws of the momentum~(\ref{scheme_momentum}) 
and the center  of mass motion~(\ref{scheme_center}).  
\end{theorem}

\begin{remark} 
The conservation of mass~(\ref{equation_discrete_mass}) provides the 
the consistency condition for the relations  
\begin{equation}         \label{scheme_mass} 
{  \hat{r}  ^{n+1}  \over  n+1   }  
- {     r   ^{n+1}  \over  n+1   }  
=  \tau _j R u ^{(0.5) } , 
\qquad 
{  r _{+}  ^{n+1} \over  n+1   }   
- {  r   ^{n+1}  \over  n+1   }  
=  { h_i  \over \rho }     .  
\end{equation}
The first relation follows from the equation~(\ref{equation_discrete_coordinate}) 
and the definition  of $R$. 
The second relation approximates $ ds / \rho  = r^n d r $, 
the differential   analog of~(\ref{Lagrange_coordinate}).   
From the consistency of the relations~(\ref{scheme_mass})  
it follows that  $   r   ^{n+1} $ plays a role of a discrete potential 
for the equation~(\ref{equation_discrete_mass}).  
Note that the mesh for the  time  and the mass Lagrange coordinate can be nonuniform. 
\end{remark}


For $ \gamma =  { 1 + 2 /  d } $  
we are interested to preserve the additional conservation laws~(\ref{Lagrange_additional_1}) 
and~(\ref{Lagrange_additional_2}).     
We can use the remaining freedom to chose the discrete equation of state. 
Let us look for an equation of state which gives us the following difference analog 
of the additional conservation law~(\ref{Lagrange_additional_1}): 
\begin{equation}          \label{scheme_additional_1} 
\left[ 
2t  \left( \varepsilon   + { < u^2 >     \over 2}       \right)   
-  < r u >     
\right]    _t  
+ 
  [   R   p  _{*} ^{(\alpha)}     ( 2t ^{(0.5)}    u  ^{(0.5)}     -    r ^{(0.5)}  )     ]  _s  
=   0    . 
\end{equation}
It leads to  the following specific  internal energy equation
\begin{equation}      \label{energy_definition} 
\varepsilon   ^{(0.5)}  = 
{  p ^{(\alpha)}   \over   \gamma -1  }   \left(   { 1\over \rho}   \right)  ^{(0.5)}
- {  \tau ^2  \over 8}       <  (  u_t  ) ^2 >  
+ { 1 \over 2 } p ^{(\alpha)}   
\left[   r ^{(0.5)}   R   -     (   r  ^{n+1}  ) ^{(0.5)}   \right]  _s    . 
\end{equation} 
We will take it as the {\it discrete equation of state}, which approximates~(\ref{ideal_gas}).

In this case we also get a difference analog 
of the second additional conservation law~(\ref{Lagrange_additional_2}) as 
\begin{multline}        \label{scheme_additional_2} 
\left[ 
t^2   \left( \varepsilon   + {  < u^2 >     \over 2}       \right)   
-   t  < r u >       
+  {  < r^2 >  \over 2 }    
+  { \tau ^2 \over 8}      { < u^2 > }
\right]    _t  
\\
\\
+ 
  [     R      p  _{*} ^{(\alpha)}   
(  (  t^2 )  ^{(0.5)}       u  ^{(0.5)}     -  t   ^{(0.5)}    r    ^{(0.5)}    )     ]  _s  
=   0   . 
\end{multline} 
Note that it has a correcting term  $ { \tau ^2  \over 8}     { < u^2 >  }   $, 
which disappears in the continuous limit.

Thus, we obtained the difference scheme~(\ref{equation_preliminary})  
supplemented by the discrete state equation~(\ref{energy_definition}). 
In this scheme the pressure values  $p$ and $ \hat{p}$  
appear only as a weighted value $  p^{(\alpha)} $, 
i.e. $ \alpha$  has no longer meaning of a parameter.   
We can consider this value 
as the  pressure in the midpoint of the cell $ (  t  _{j+1/2}  ,   s  _{i+1/2} )  $, 
i.e.  for  $\alpha = 0.5$.

Thus, we arrive at the following result:

\begin{theorem}    \label{second_theorem} 
The scheme~(\ref{equation_preliminary}) with the discrete equation of state~(\ref{energy_definition}) 
possesses the properties given in Theorem~\ref{first_theorem}  
as well as the additional conservation laws~(\ref{scheme_additional_1}) and~(\ref{scheme_additional_2})
for $ \gamma =  { 1 + 2 /  d } $. 
\end{theorem}


With the help of 
\begin{equation*}   
\left[   r ^{(0.5)}   R   -     (   r   ^{n+1}  )  ^{(0.5)}    \right] 
= 
\left\{ 
\begin{array}{ll} 
0 , &   n  =  0 ; \\
 &  \\
{\displaystyle     - { 1\over 4}  ( \hat{r} - r ) ^2   
=    - { 1\over 4}  \tau ^2   ( u ^ {(0.5)}  ) ^2  }   , &   n  =  1 ; \\ 
 &  \\
{\displaystyle    - {   1  \over 3}   (   \hat{r} + r   )   ( \hat{r} - r ) ^2  
=   - {   2    \over 3}    \tau ^2  r ^ {(0.5)}     ( u ^ {(0.5)}  ) ^2  }   , &   n  =  2 \\ 
\end{array} 
\right. 
\end{equation*}   
we conclude that the additional terms in the discrete equation of state~(\ref{energy_definition}) 
represent    $ O ( \tau ^2   )  $ correction to the continuous equation of state~(\ref{ideal_gas}).

\section{Concluding remarks} 

\label{conclusion}

In the present paper we considered one-dimensional flows of a polytropic gas. 
There were derived difference schemes 
which in addition to conservation laws of mass and energy 
(as well as conservation of momentum and the center of mass motion for the plain one-dimensional flows) 
posses two additional conservation laws for the special values of the adiabatic exponent $\gamma = 1 + 1/d $.

The case of the plain one-dimensional flows was considered  in~\cite{Korobitsyn}, 
where the mesh for  the mass Lagrangian coordinate  was uniform. 
Here the results of~\cite{Korobitsyn}  are extended to nonuniform spatial meshes.

\section*{Acknowledgements}

The research was supported by Russian Science Foundation Grant no. 18-11-00238 
"Hydrodynamics-type equations: symmetries, conservation laws, invariant difference schemes".



\begin{thebibliography}{9999}







\par
\bibitem{Chernyi} 
G.~G.~Chernyi  (1988)
{\it Gas dynamics}, 
Nauka, Moscow  (in Russian).



\par
\bibitem{Chorin} 
A.~J.~Chorin  and J.~E.~Marsden 
(1990) 
{\it  A Mathematical Introduction to Fluid Mechanics}, 
Springer-Verlag. 

 
\par
\bibitem{Furihata}  
D.~Furihata and T.~Matsuo 
(2011) 
{\it Discrete Variational Derivative Method: 
A Structure-Preserving Numerical Method for Partial Differential Equations}, 
CRC Press.





\par
\bibitem{Hairer}  
E.~Hairer, C.~Lubich and G.~Wanner 
(2006)
{\it Geometric Numerical Integration: 
Structure-Preserving Algorithms for Ordinary Differential Equations}, 
Springer.


\par
\bibitem{Ibragimov2}
N.~H.~Ibragimov
(1973)
Conservation laws in hydrodynamics,  
{\it Dokl. Akad. Nauk SSSR} 
{\bf 210}  (6)  (1973) 1307--1309. 



\par
\bibitem{Ibragimov}
N.~H.~Ibragimov 
(1985) 
{\it Transformation Groups Applied to Mathematical Physics}, 
Reidel, Boston. 






\par
\bibitem{Korobitsyn}
V.~A.~Korobitsyn    
(1989)  
Thermodynamically matched difference schemes
{\it Zh. Vychisl. Mat. Mat. Fiz.}  {\bf 29} (2)    309--312.  





\par
\bibitem{Landau}
L.~D.~Landau and E.~M.~Lifshitz  
(1987) 
{\it Fluid Mechanics},  2nd. ed., 
Pergamon Press. 



\par
\bibitem{LeVeque}  
R.~J.~LeVeque  
(1992) 
{\it  Numerical Methods for Conservation Laws}, 
Birkhauser-Verlag.




\par
\bibitem{Laney} 
C.~B.~Laney    
(1998)  
{\it    Computational {Gas}dynamics},   Cambridge Univ. Press.



\par
\bibitem{Moskalkov1}  
M.~N.~Moskal'kov  
(1980) 
A completely conservative scheme of gas dynamics, 
{\it Zh. Vychisl. Mat. Mat. Fiz.}  {\bf 20} (1),  162--170; 
{\it U.S.S.R. Comput. Math. Math. Phys.}  {\bf  20} (1), 177--187.  


\par
\bibitem{Moskalkov2} 
M.~N.~Moskal'kov 
(1981)
Using the dispersion properties of a difference scheme of gas dynamics, 
{\it Zh. Vychisl. Mat. Mat. Fiz.}   {\bf  21} (5),  1257--1263; 
{\it U.S.S.R. Comput. Math. Math. Phys.} {\bf  21} (5), 182--188.  




\par
\bibitem{Ovsiannikov} 
L.~V.~Ovsiannikov  
(2003)  
{\it Lectures on the gas dynamics equations}, 
Institute of computer studies, Moscow--Izhevsk    (in Russian). 




\par
\bibitem{Popov} 
Yu.~P.~Popov and A.~A.~Samarskii  
(1969)  
Completely conservative difference schemes 
{\it Zh. Vychisl. Mat. Mat. Fiz.} {\bf  9}  (4)  953--958. 



 
 
 





\par
\bibitem{Samarskii} 
A.~A.~Samarskii and Yu.~P.~Popov  
(1980) 
{\it   Difference methods for solving problems of gas dynamics},  
Nauka, Moscow     (in Russian).



\par
\bibitem{Shashkov}  
M.~Shashkov  
(1996) 
{\it Conservative Finite-Difference Methods on General Grids}, 
CRC Press, Boca Raton, Fl.






\par
\bibitem{Terentev} 
E.~D.~Terentev and Yu.~D.~Shmyglevskii 
(1975) 
A complete system of equations in divergence form for the dynamics of an ideal gas 
{\it  Zh. Vychisl. Mat. Mat. Fiz.} 
{\bf 15} (6)   1535--1544. 




\par
\bibitem{Toro} 
E.~F.~Toro   
(1997) 
{\it Riemann Solvers and Numerical Methods for Fluid Dynamics}, 
Springer-Verlag, Berlin-Heidelberg. 




\end{thebibliography}
\end{document}